\documentclass[12pt,draft]{article}
\usepackage{amssymb, amsmath, amsthm, verbatim, cite, lscape}
\usepackage[cp1251]{inputenc}
\usepackage[english]{babel}

\newtheorem{theorem}{Theorem}
\newtheorem{corol}{Corollary}

\newcounter{parag}
\newcommand{\sect}[1]
{\refstepcounter{parag}
\begin{center} { \bf\theparag. #1} \end{center}}

\newtheorem{lemma}{Lemma}[parag]

\theoremstyle{definition}
\newtheorem*{prf}{Proof}
\newtheorem*{prf1}{Proof of Theorem 1}

\newtheorem*{prf3}{Proof of Theorem 2}
\newtheorem*{prf4}{Proof of Theorem 3}

\theoremstyle{definition}

\begin{document}
\title{The time complexity of some algorithms for generating the spectra of finite simple groups}
\author{Alexander Buturlakin\thanks{buturlakin@math.nsc.ru; Sobolev Institute of Mathematics, 4~Acad. Koptyug avenue, 630090 Novosibirsk, Russia.}}
\date{}
\maketitle
\abstract{The spectrum $\omega(G)$ is the set of orders of elements of $G$. We consider the problem of generating the spectrum of a finite nonabelian simple group $G$ given by the degree of $G$ if $G$ is an alternating group, or the Lie type, Lie rank and order of the underlying field if $G$ is a group of Lie type.}



\sect{Introduction}

The spectrum $\omega(G)$ is the set of orders of elements of a group $G$. Let $G$ be a finite nonabelian simple group. We discuss whether there exists an algorithm that given the degree of $G$ if it is an alternating group, or the Lie type, Lie rank and order of the underlying field if $G$ is a group of Lie type, outputs the spectrum of $G$. It is known that the size of the spectrum of an alternating group is an exponential function of its degree (see Lemma~\ref{l:ErdTur} below). The size of the spectrum of a group of Lie type happens to be an exponential function of its Lie rank (Lemma~\ref{l:MuLieMuSn} below) and is a polynomial function of the order of the field (this immediately follows from the description of the spectra of these groups). So we are looking for an algorithm with running time bounded by a polynomial in the length of the output. Define the length of a set of positive integers $\mathcal{S}$ to be the sum of $\log s$ (note that $\log$ will always mean the natural logarithm) where $s$ runs over~$\mathcal{S}$.

Let $\mu(G)$ be the set of elements of $\omega(G)$ that are maximal with respect to divisibility. Since the spectrum is closed under taking divisors, it is uniquely determined by any subset of $\omega(G)$ that contains~$\mu(G)$.

\begin{theorem}\label{t:Alt} Let $G$ be a symmetric or alternating group of degree~$n$. There is an algorithm that given $n$ outputs $\mu(G)$ and $\omega(G)$ in time bounded by a polynomial of the length of the output.
\end{theorem}

Thus, the required algorithm exists in the case of alternating groups. As for the groups of Lie type, we have a weaker result. For the exceptional groups, we still have a polynomial algorithm but only for~$\mu(G)$.

\begin{theorem}\label{t:LieExcep} Let $G$ be a finite exceptional group of Lie type over a field of order~$q$. There is an algorithm that given the Lie type of $G$ and $q$ outputs $\mu(G)$ in time bounded by a polynomial of the length of the output.
\end{theorem}

Theorem~\ref{t:LieExcep} directly follows from the existing description of the spectra of the exceptional groups. For more details see the proof of this theorem in Section~3.

In the case of classical groups, we only prove the existence of a quasi-polynomial algorithm.

\begin{theorem}\label{t:Lie} Let $G$ be a finite classical group of Lie rank $n$ over a field of order~$q$. Let $m$ be the length of $\mu(G)$. There is an algorithm that given the Lie type of $G$, $n$ and $q$ outputs $\mu(G)$ in time $m^{O(\sqrt{\log \log m})}$.
\end{theorem}

Observe that we actually do not provide any new algorithm for the spectrum generation. We just analyze the efficiency of the most straightforward algorithms provided by the known description of spectra of groups of Lie type (see \cite{08But.t, 10But.t} for the classical groups, and \cite{18But.t} and references therein for the exceptional groups) and the description of conjugacy classes of alternating groups. We also want to note that we have no evidence for the existence or non-existence of a polynomial-time algorithm for the classical groups.

The results of this paper complement the main theorem of \cite{19ButVas}, which is concerned with the following problem. Given a finite set $\mathcal M$ of natural numbers, determine whether there exists a finite simple group $G$ whose spectrum coincide with the set of divisors of elements of~$\mathcal M$. It is proved in \cite{19ButVas} that there is a polynomial-time algorithm which returns a unique candidate for $G$ (more precisely, its parameters in the sense described above) or says that there is no such group. This candidate $G$ satisfies the following conditions: $\mathcal M$ is a subset of $\omega(G)$ and if $H$ is a finite simple group whose spectrum differs from the spectrum of $G$, then the set of divisors of elements of $\mathcal M$ is not the spectrum of~$H$. Thus, to complete the task it remains to verify whether $\omega(G)$ is equal to the set of divisors of elements of $\mathcal M$, or equivalently, whether $\mu(G)$ is a subset of $\mathcal M$. Hence the following statements are corollaries of the main result of \cite{19ButVas} and Theorems~\ref{t:Alt}, \ref{t:LieExcep} and \ref{t:Lie}.

\begin{corol}\label{c:AlgAlt} There is an algorithm that given a set $\mathcal M$ of positive integers  outputs the degree of an alternating group whose spectrum coincides with the set of divisors of elements of $\mathcal M$, or says that there is no such a group. The running time of the algorithm is polynomial in the length of $\mathcal M$.
\end{corol}

\begin{corol}\label{c:AlgLieExcep} There is an algorithm that given a set $\mathcal M$ of positive integers  outputs a finite simple exceptional group of Lie type (the Lie type, Lie rank and order of the underlying field) whose spectrum coincides with the set of divisors of elements of $\mathcal M$, or says that there is no such group. The running time of the algorithm is polynomial in the length of $\mathcal M$.
\end{corol}

\begin{corol}\label{c:AlgLie} There is an algorithm that given a set $\mathcal M$ of positive integers  outputs a finite simple classical group of Lie type (the Lie type, Lie rank and order of the underlying field) whose spectrum coincides with the set of divisors of elements of $\mathcal M$, or says that there is no such a group. The running time of the algorithm is $m^{O(\sqrt{\log \log m})}$, where $m$ is the length of $\mathcal M$.
\end{corol}

\sect{Alternating groups}

Denote by $Sym_n$ and $Alt_n$ the symmetric and alternating group of degree~$n$. We write $g(n)$ for Landau's function of a positive integer $n$, that is, the largest order of an element of the group $Sym_n$.

\begin{lemma}\label{l:landsu} If $n\geqslant 810$, then $0.99\leqslant\log g(n)/\sqrt{n\log n}\leqslant 1.08$.
\end{lemma}

\begin{prf} Directly follows from the bounds for $g(n)$ introduced in \cite[Theorems 1, 2]{89MasNicRob}.
\end{prf}

For a positive integer $n$, denote by $\sigma(n)$ the divisor function of $n$, that is, the number of divisors of~$n$. The following is a direct corollary of \cite[Theoreme 1]{83NicRob}.

\begin{lemma}\label{l:BoundForDivisors} Put $$n_0=2^5\cdot 3^3\cdot5^2\cdot7\cdot11\cdot13\cdot 17\cdot19$$ and $$\alpha_0=\frac{\log\sigma(n_0)\log\log n_0}{\log n_0}.$$ If $n\geqslant 2$, then $$\log \sigma(n)\leqslant \alpha_0\frac{\log n}{\log\log n}.$$
\end{lemma}

\begin{lemma}\label{l:ErdTur}

\begin{itemize}
\item[$(1)$] $|\omega(Sym_n)|=\exp\left(\frac{2\pi}{\sqrt{6}}\sqrt{\frac{n}{\log n}}+O\left(\frac{\sqrt{n}\log\log n}{\log n}\right)\right).$

\item[$(2)$] $|\mu(Sym_n)|\geqslant\exp\left(c\sqrt{\frac{n}{\log n}}+o\left(\sqrt{\frac{n}{\log n}}\right)\right)$ for a positive $c$ greater than $\left(\frac{2\pi}{\sqrt{6}}-2.16\alpha_0\right)>0.26$, where $\alpha_0$ is defined in Lemma~\ref{l:BoundForDivisors}.
\end{itemize}
\end{lemma}

\begin{prf} Statement $(1)$ is proved in {\cite[Theorem I]{68ErdTur}}. To prove the second statement, let us bound the number $\sigma(m)$ of divisors  of an element $m$ of $\mu(Sym_n)$. Since the bound in Lemma~\ref{l:BoundForDivisors} is an increasing function, we can write $$\log\sigma(m)\leqslant \alpha_0\frac{\log g(n)}{\log\log g(n)}.$$ It follows from Lemma~\ref{l:landsu}, that $$\log\sigma(m)\leqslant \alpha_0\frac{1.08\sqrt{n\log n}}{\log (1.08\sqrt{n\log n})}= \frac{2\cdot1.08\alpha_0\sqrt{n\log n}}{2\log 1.08+\log n+\log\log n}\leqslant 2.16\alpha_0 \sqrt{\frac{n}{\log n}}.$$  Finally, the inequalities $$|\omega(Sym_n)|\geqslant|\mu(Sym_n)|\geqslant\frac{|\omega(Sym_n)|}{\max_{m\in\mu(Sym_n)}\sigma(m)}$$ give the required asymptotic.
\end{prf}

Given a positive integer $n$, define $l(n)$ to be the minimal degree such that $Sym_{l(n)}$ contains an element of order $n$. That is, $l(n)=0$ if $n=1$, and $l(n)=p_1^{\alpha_1}+\dots+p_k^{\alpha_k}$ if $n\neq 1$ and $n=p_1^{\alpha_1}\dots p_k^{\alpha_k}$, where $p_1$, $\dots$, $p_k$ are pairwise distinct primes and every $\alpha_i>0$.

\begin{prf1} Let us first prove the statement for the whole spectrum in the case of symmetric group. The algorithm proceeds as follows. First, it generates the list of all primes up to $n$ using the sieve of Eratosthenes. Denote them by $p_1<p_2<\dots <p_k$. Put $P_i=\{p_i^a|\ a\geqslant 0\text{ is an integer and } p_i^a\leqslant n\}$. Next the algorithm generates the sets $P_i$ for $1\leqslant i\leqslant k$. The elements of $\omega(Sym_n)$ are parameterized by elements of the direct product of the sets~$P_i$. Now the algorithm picks some $x_1$ from $P_1$. If $n-l(x_1)<p_2$, then the algorithm outputs $(x_1, 1,\dots, 1)$. Otherwise it proceeds to $P_2$. Assume that $x_1$, $\dots$, $x_i$ are chosen. If $i=k$ or $n-l(x_1\dots x_i)<p_{i+1}$, then the algorithm outputs $(x_1,\dots, x_i, 1,\dots, 1)$. Otherwise it chooses $x_{i+1}$ from  $P_{i+1}$ such that $l(x_1\dots x_{i+1})\leqslant n$.

Obviously this algorithm visits every element of $\omega(Sym_n)$ once and the number of steps required to write down such an element is at most $k$, which is less than $n$. Hence, it remains to show that the time required to find all primes not exceeding $n$ is polynomial in output.

By Lemma~\ref{l:ErdTur}, there exists a number $N$ such that for every $n>N$ we have $$\exp\left(\frac{3\pi}{\sqrt{6}}\sqrt{\frac{n}{\log n}}\right)\geqslant|\omega(Sym_n)|\geqslant \exp\left(\frac{\pi}{\sqrt{6}}\sqrt{\frac{n}{\log n}}\right).$$  Hence the length of the output is an exponential function of $n$. The sieve of Eratosthenes requires $O(n\log\log n)$ steps, and the theorem is proved for symmetric groups.

To obtain an analogous algorithm for alternating group, one should change the definition of
$P_1$: now $P_1$ = $\{2^a+2| a\geqslant0 \text{ is an integer and } 2^a+2\leqslant n\}$. The rest of the algorithm
remains unchanged.

Finally, Lemma~\ref{l:ErdTur} implies that there exist polynomials $p_1$ and $p_2$ such that $$|\omega(Sym_n)|\leqslant p_1(|\mu(Sym_n)|) \text{ and } |\omega(Alt_n)|\leqslant p_2(|\mu(Alt_n)|).$$ Hence the elements of $\mu(G)$ where $G$ is $Alt_n$ or $Sym_n$ can be constructed from $\omega(G)$
in time polynomial in the length of~$\mu(G)$. This completes the proof.

\end{prf1}

\sect{Groups of Lie type}

\begin{prf3} The spectra of all finite simple exceptional groups of Lie type are described (see, for example, \cite{18But.t} and references therein). It follows from these descriptions that there exists a constant $C$ not depending on $G$ such that the spectrum of $G$ contains a subset $\nu(G)$ with the following properties: $\mu(G)$ is a subset of $\nu(G)$, the cardinality of $\nu(G)$ is at most $C$ and every element of $\nu(G)$ can be computed in time polynomial in $\log q$, where $q$ is the order of the underlying field of~$G$. Hence $\nu(G)$ can be computed in time bounded by a polynomial of $\log q$ and so can be $\mu(G)$. This proves the theorem.
\end{prf3}

Before proving Theorem~\ref{t:Lie}, we need to introduce some notations and preliminary results.

Fix an integer $a$ with $|a|>1$. If $s$ is an odd integer coprime to $a$, define $e(s,a)$ to be the multiplicative order of $s$ modulo $a$. If $a$ is odd, put $e(2,a)=1$ if $a\equiv 1$ modulo $4$, and $e(2,a)=2$ otherwise. A prime $r$ is said to be a primitive prime divisor of
$a^i-1$ if $e(r,a)=i$. We write $r_i(a)$ to denote some primitive prime divisor of $a^i-1$, if such
a prime exists, and $R_i(a)$ to denote the set of all such divisors. In \cite{86Ban, Zs}, it is proved
that primitive prime divisors exist for almost all pairs $(a,i)$.

\begin{lemma}\label{l:zsigmondy} Let $a$ be an integer and $|a| > 1$. For every positive integer $i$ the set $R_i(a)$ is nonempty, except for the pairs $(a, i)\in \{(2, 1), (2, 6), (-2, 2), (-2, 3), (3, 1), (-3, 2)\}$.
\end{lemma}




For a pair of integers $a$ and $b$, denote by $(a,b)$ and $[a,b]$ their greatest common divisor and least common multiple. If $b$ is a prime, then denote by $(a)_b$ the $b$-part of $a$, that is, the greatest power of $b$ that divides~$a$, and by $(a)_{b'}$ the ratio $|a|/(a)_b$.

The following lemma is elementary.

\begin{lemma}\label{l:gcd} Let $a$ be an integer such that $|a|>1$ and $s$, $t$ positive integers.

$1)$ $(a^s-1, a^t-1)=a^{(s,t)}-1;$

$2)$ $(a^s+1, a^t-1)=\begin{cases}a^{(s,t)}+1 \text{ if } \frac{s}{(s,t)} \text{ is odd and } \frac{t}{(s,t)} \text{ is even},\\ (2,a-1) \text{ otherwise}.\end{cases}$
\end{lemma}

If $r_i(a)$ exists, let $r^*_i(a)$ be the greatest power of $r_i(a)$ dividing $a^i-1$. Also define $r_3^*(-2)=r_6^*(2)=9$ and $r_2^*(-3)=8$. The following statement is a direct consequence of Lemma~\ref{l:gcd}.

\begin{corol}\label{c:DivOfLcm} Let $s$ be a prime, $\alpha$ a positive integer, and $a$ an integer such that $|a|>1$. Suppose that $r^*_{s^\alpha}(a)$ is defined. If $r^*_{s^\alpha}(a)$
divides $$A=[a^{n_1}-\varepsilon_1, \dots, a^{n_l}-\varepsilon_l],$$ where $\varepsilon_i\in\{ 1, -1\}$, then one of the following statements holds:

$1)$ there exists $i\in\{ 1,\dots,l\}$ such that $s^\alpha$ divides $n_i$;

$2)$ $s=2$ and there exists $i\in\{ 1,\dots,l\}$ such that $\varepsilon_i=-1$ and $s^{\alpha-1}$ divides $n_i$.

If $(r^*_{s^\alpha}(a), a-1)\neq 1$, then either $s^\alpha=2$ and $r^*_{s^\alpha}(a)$ is a power of $2$, or $s^\alpha=3$ and $r^*_{s^\alpha}(a)$ is a power of $3$.
\end{corol}


Denote by $p(n)$ the number of partitions of $n$. For a group of Lie type over a field of characteristic $p$, denote by $\omega_{p'}(G)$ the set of orders of semisimple elements, that is, elements whose order is coprime to $p$. Put $\mu_{p'}(G)=\mu(\omega_{p'}(G))$. Following \cite{15Vas}, we use the notation $\operatorname{prk}(G)$ for the dimension of $G$ if $G$ is a linear or unitary group and the Lie rank of $G$ if $G$ is an orthogonal or symplectic group.

\begin{lemma}\label{l:MuLieMuSn} Let $G$ be a simple classical group with $\operatorname{prk}(G)=n$. Then $$|\mu_{p'}(G)|\geqslant|\mu(Sym_n)|$$ unless $G=PSU_n(2)$  and always  $$|\mu(G)|\geqslant|\mu(Sym_n)|.$$
\end{lemma}

\begin{prf}
We will prove the lemma by constructing an injection $\phi$ from $\mu(Sym_n)$ to~$\mu_{p'}(G)$ for groups distinct from $PSU_n(2)$ and to $\mu(G)$ otherwise.

Assume that $G=PSL_n^\varepsilon(q)$, where $\varepsilon\in\{ +,- \}$. Suppose that $n_1,\dots, n_s$ are positive integers such that $n_1+\dots+n_t\leqslant n$. Define
$$c(n_1,\dots,n_t)=
\begin{cases}
1& \text{if } t+n-(n_1+\dots+n_t)>2;\\
(n/(n_1,n_2), \varepsilon q-1) & \text{if } t=2 \text{ and }n=n_1+\dots+n_t;\\
(n,\varepsilon q-1)& \text{if } t=1 \text{ and } n=1+n_1+\dots+n_t;\\
(q-\varepsilon1)((n,\varepsilon q-1)& \text{if } t=1 \text{ and } n=n_1+\dots+n_t.\\
\end{cases}$$
In any case, $\pi(c(n_1,\dots,n_t))\subseteq\pi(\varepsilon q-1)$.

According to \cite[Corollary 3]{08But.t},
for every $n_1,\dots,n_t$ such that $n_1+\dots+n_t\leqslant n$,
$\omega(G)$ contains the number $$[(\varepsilon q)^{n_1}-1,\dots, (\varepsilon q)^{n_t}-1]/c(n_1,\dots, n_t),$$
and every element of $\mu_{p'}(G)$ has such a form.

Let $$a=p_1^{\alpha_1}\dots p_t^{\alpha_t}$$ for prime numbers  $p_1<\dots<p_t$  be an element of $\mu(Sym_n)$. Since $l(a)\leqslant n$, the number of the form $$b_a=[(\varepsilon q)^{p_1^{\alpha_1}}-1, \dots, (\varepsilon q)^{p_t^{\alpha_t}}-1]/c(p_1^{\alpha_1},\dots, p_t^{\alpha_t})$$ is an element of $\omega(G)$.

Assume that for some $i$, the number $r^*=r^*_{p_i^{\alpha_i}}(\varepsilon q)$ is defined but does not divide~$b_a$.  Then
$r^*$ is not coprime to $c(p_1^{\alpha_1},\dots, p_t^{\alpha_t})$, in particular, $t\leqslant 2$. By Corollary \ref{c:DivOfLcm}, it follows that $p_i^{\alpha_i}=2$ or $3$ and $r^*$  is a power of $p_i^{\alpha_i}$. If $t=2$, then  $n=p_i^{\alpha_i}+p_j^{\alpha_j}$ is coprime to $p_i^{\alpha_i}$, and so $(c(p_i^{\alpha_i}, p_j^{\alpha_j}), r^*) =1$. If $t=1$, then $n=2$ or $3$ (the lemma in this cases can be easily checked and we will assume in what follows that $n>3$). So we may assume that every $r^*_{p_i^{\alpha_i}}(\varepsilon q)$ divides~$b_a$ if exists.

Suppose that $r^*_{p_i^{\alpha_i}}(\varepsilon q)$ is defined for all $i$ and let $A$ be an element of $\mu_{p'}(G)$ divisible by $b_a$. We define $\phi(a)=A$.

Suppose that there exists an element $a'=r_1^{\beta_1}\dots r_u^{\beta_u}$ of $\mu(Sym_n)$ such that $a'\neq a$, all
$r^*_{r_j^{\beta_j}}(\varepsilon q)$ are defined, and $\phi(a)=\phi(a')$. Let $$A=[(\varepsilon q)^{m_1}-1,\dots, (\varepsilon q)^{m_s}-1]/c(m_1,\dots, m_s)$$ with $m_1+\dots+m_s\leqslant n$.  It follows from Corollary~\ref{c:DivOfLcm} that for each $p_i^{\alpha_i}$, there exists a number $m_{k}$ divisible by $p_i^{\alpha_i}$. The same is true for each $r_j^{\beta_j}$. Since  $m_1+\dots+m_s\leqslant n$, the least common multiple of $m_1$, $\dots$, $m_s$ lies in $\mu(Sym_n)$. By maximality of $a$ and $a'$, this number is equal to both $a$ and $a'$, which is a contradiction.



Assume that some $r^*_{p_i^{\alpha_i}}(\varepsilon q)$ does not exists. By Lemma~\ref{l:zsigmondy}, we have $G=PSU_n(2)$ and $p_i^{\alpha_i}=2$ (in fact, $p_1^{\alpha_1}=2$). Observe that the cases when $n<8$ can be checked by direct computations. Since $2a$ is not an element of $\mu(Sym_n)$, the sum of $p_i^{\alpha_i}$ for $i>1$ is $n-3$ or $n-2$. Define $\phi(a)$ to be $$2\left[2^{p_2^{\alpha_2}}+1, \dots, 2^{p_t^{\alpha_t}}+1\right]/d \text{ if this sum is } n-2,$$ and as $$4\left[2^{p_2^{\alpha_2}}+1, \dots, 2^{p_t^{\alpha_t}}+1\right]/d \text{ otherwise},$$ where in both cases $d=1$ if $t>1$, and $d=(n,3)$ otherwise.   Using \cite[Corollary 3]{08But.t} and the fact that $r^*_{p_i^{\alpha_i}}(\varepsilon q)$ divides $\phi(a)$ for $i>1$, it is easy to see that $\phi(a)\in\mu(G)$ and the resulting
$\phi$ is an injection from $\mu(Sym_n)$ to $\mu(PSU_n(2))$.
Thus the lemma is proved for all linear and unitary groups.

Assume that $G=PSp_{2n}(q)$ or $\Omega_{2n+1}(q)$. Again let $a=p_1^{\alpha_1}\dots p_t^{\alpha_t},$ where $p_1$, $\dots$, $p_t$ are pairwise distinct primes,  be an element of $\mu(Sym_n)$.
By \cite[Corollaries 2, 3, 6]{10But.t},  $\omega(G)$ contains the number
\begin{equation}\label{eq2'}
b_a=[q^{p_1^{\alpha_1}}+(-1)^{p_1}, \dots, q^{p_t^{\alpha_t}}+(-1)^{p_t}]/(2,q-1).\end{equation} Observe that $b_a$ is divisible by all $r^*_{p_i^{\alpha_i}}(q)$ for $p_i\neq 2$ and also by
$r^*_{2p_i^{\alpha_i}}(q)$ if $p_i=2$. Define $\phi(a)$ to be an element of $\mu_{p'}(G)$ divisible by $b_a$. It follows from \cite[Corollaries 2, 3, 6]{10But.t} that $\phi(a)$ has the form \begin{equation}\label{eq2}[ q^{m_1}-\varepsilon_1, \dots, q^{m_s}-\varepsilon_s ]/c,\end{equation} where $\varepsilon_i\in\{1, -1\}$ and $c=1,2$. It follows from Corollary~\ref{c:DivOfLcm}, that for each $i$ there exists index $j$ such that $p_i^{\alpha_i}$ divides $m_j$. As before, the assumption that $\phi$ is not injective leads to a contradiction with the maximality of $a$ in $\mu(Sym_n)$.

Let $G=P\Omega_{2n}^\varepsilon(q)$, where $\varepsilon\in\{+, -\}$. By \cite[Corollaries 8, 9]{10But.t}, the number of the form (\ref{eq2'}), but with possibly $(4,q-\varepsilon)$ instead of $(2,q-1)$ in the denominator, lies in $\omega(G)$ if and only if either $l(a)<n$, or the number of pluses among $(-1)^{p_i}$ is even when $\varepsilon=+$ and odd otherwise. Hence we should worry about the signs in the definition of $b_a$ in this case. If $l(a)<n$, then the previous definition of $b_a$ works.
Suppose that $l(a)=n$. If $t\leqslant 2$ or $p_1\neq 2$, then one can substitute $q^{p_i^{\alpha_i}}-1$ by $q^{p_i^{\alpha_i}}+1$ for one odd prime $p_i$ if necessary. If $n=a$ is a power of $2$,
then $b_a=(q^n+1)/(2,q-1)$ if $\varepsilon=-$ and
$b_a=(q^n-1)_{2'}$ if $\varepsilon=+$. Again $b_a$ is divisible by $r^*_{p_i^{\alpha_i}}(q)$ or $r^*_{2p_i^{\alpha_i}}(q)$ for each odd $p_i$ and also by $r^*_{2p_i^{\alpha_i}}(q)$ if $p_i=2$
and $(n,a,\varepsilon)\neq (2^l,2^l,+)$. As before, define $\phi(a)$ to be an element of $\mu_{p'}(G)$ divisible by $b_a$.

To see that $\phi$ is bijection, it is sufficient to show that $\phi(a)\neq \phi(n)$ when $\varepsilon=+$, $n$ is a power of $2$ and $a\neq n$. Suppose that $\phi(a)=\phi(n)=A$. Then $A$ has the form as in \eqref{eq2} but with $c=1,2,4$. Reasoning as above, we conclude that $a=[m_1,\dots,m_s]$.
Since $A$ is divisible by $r_n^*(q)$ then there exists index $j$ such that $n/2$ divides $m_j$. If $m_j=n$, then $a=n$, so $m_j=n/2$. Since $A$ is divisible by $r_{n/2}^*(q)$ then there exists index $j_1$ such that $n/4$ divides $m_{j_1}$. Repeating this
argument until $n/2^l\geqslant 4$, we see that each of $m_1,\dots,m_s$ divides $m_j$ and so $a=n/2$, which contradicts maximality of $a$.  The lemma is proved.
\end{prf}

\begin{prf4}  Let $G$ be a classical group with $\operatorname{prk}(G)=n$ over a field of characteristic $p$ and order $q$. By Lemma~\ref{l:MuLieMuSn}, we have $|\mu(G)|\geqslant|\mu(Sym_n)|$. Hence there exist $\varepsilon>0$ and $M$ such that for all $n>M$, we have \begin{equation}\label{eq3}|\mu(G)|>\exp\left(\varepsilon\sqrt{\frac{n}{\log n}}\right).\end{equation}

If $G$ is a linear or unitary group, then by \cite[Corollary 3]{08But.t} there is a subset $\nu(G)$ of $\omega(G)$ containing $\mu(G)$ whose elements are parameterized by partitions of $n$ of the form \begin{equation}\label{eq4}n=(p^{\alpha}-1)+n_1+\dots+n_s,\end{equation} where first term can be not presented.
There are at most $1+\log_2 (n+1)$ choices for the first term (including its absence), and the number of choices for $n_1$, $\dots$, $n_s$ obviously not exceeds $p(n)$ for each choice of the first term. Therefore we have $$|\nu(G)|<p(n)\log_2 (2n+2).$$ Every element of $\nu(G)$ can be computed in time polynomial in $n\log q$. Hence the bound (\ref{eq3}) provides the required bound on the running time of the algorithm that just generates all partitions of the form (\ref{eq4}) and calculates corresponding elements of~$\nu(G)$.

Let $G$ be a symplectic or orthogonal group. By \cite[Corollaries 2, 3, 6]{10But.t}, there is a subset $\nu(G)$ of $\omega(G)$ containing $\mu(G)$ whose elements are parameterized by partitions of $n$ of the form $$n=f(p)+\overline m_1+\overline m_2,$$ where $\overline k$ denote a partition of a number $k$ and first term $f(p)$ can be not presented. The function $f(p)$ depends on the Lie type of $G$ and its characteristic and is equal to  $\frac{p^{k}+i}{2}$ for some $i\in\{1, 2, 3, 4\}$. As before, we obtain the required bound on the running time of the brute force algorithm in this case. This completes the proof.

\end{prf4}

\sect{Acknowledgments}

The author thanks M.A. Grechkoseeva and A.V. Vasil'ev for valuable comments.

The work is supported by the Mathematical Center in Akademgorodok under the agreement No.~075-15-2019-1613 with the Ministry of Science and Higher Education of the Russian Federation.

\end{document}